%% file: alc.tex
\documentclass[12pt]{article}
\usepackage{graphicx,latexsym,url}
\usepackage{amsmath,amsthm}
\usepackage[top=2cm, bottom=3cm, left=3cm, right=3cm]{geometry}

\newcommand{\myadd}{\mbox{\footnotesize $\,+\,$} }
\newcommand{\mA}{\mathcal{A}}
\newcommand{\mI}{\mathcal{I}}
\newcommand{\mL}{\mathcal{L}}
\newcommand{\mN}{\mathcal{N}}
\newcommand{\mP}{\mathcal{P}}
\newcommand{\mR}{\mathcal{R}}
\newcommand{\mS}{\mathcal{S}}

\newcommand{\mU}{\mathcal{U}}
\newcommand{\mZ}{o\mathcal{A}x}
\newcommand{\mO}{\mathcal{O}}
\newcommand{\moO}{o\mathcal{O}}
\newcommand{\moOo}{o\mathcal{O}o}
\newcommand{\moA}{o\mathcal{A}}
\newcommand{\mX}{\mathcal{X}}
\newcommand{\mxX}{x\mathcal{X}}
\newcommand{\mxXx}{x\mathcal{X}x}
\newcommand{\mAx}{\mathcal{A}x}
\def\myup{\mbox{$\uparrow$}}
\def\mydn{\mbox{$\downarrow$}}

\newtheorem{theorem}{Theorem}
\newtheorem*{lemma*}{Lemma}
\newtheorem*{corollary*}{Corollary}
\newtheorem*{obs*}{Observation}

\begin{document}
{\centering 

{\large\bf A Proof of the 2004 Albert-Grossman-Nowakowski-Wolfe Conjecture on Alternating Linear Clobber}

Xinyue Chen, Taylor Folkersen, Kamillah Hasham, 
Ryan B.\ Hayward, 
David Lee, Owen Randall, Luke Schultz, Emily Vandermeer

\verb+hayward@ualberta.ca+

}

{\bf Abstract.} 
Clobber is an alternate-turn two-player game introduced in 2001 by
Albert, Grossman, Nowakowski and Wolfe.
The board is a graph with each node colored 
black (here denoted with the token {\tt x}), 
white (token {\tt o}) or empty (token {\tt -}).
Player Left has black stones, player Right has white stones.
On a turn, a player takes one of their stones
that is adjacent to an opponent stone and 
clobbers the opponent's stone (replaces it with theirs).
Whoever cannot move loses.
Linear clobber is clobber played on a path, for example, one row of a go board.
In 2004 Albert et al.\ conjectured that,
for every non-empty even-length alternating-color linear clobber position
except {\tt oxoxox},
the first player has a winning strategy.
We prove their conjecture.


\begin{figure}[ht]
\hfill\includegraphics[scale=.5]{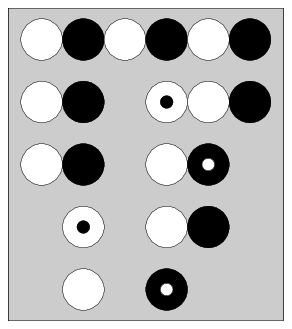}\hfill\ 
\caption{From top, a linear clobber game.}
\label{fdemo6}
\end{figure}

\section*{Introduction}
Clobber is an alternate-turn two-player combinatorial game introduced by 
Michael H.\ Albert, J.P.\ Grossman, Richard J.\ Nowakowski and David Wolfe
in 2001 \cite{AGNW}.
Left has black stones 
(here denoted with the token {\tt x}), 
Right has white stones (token {\tt o}).
Clobber can be played on any graph.
To start, each board cell (graph node) 
has three possible states (colors):
no stone, a black stone, or a white stone.
On a turn, a player selects one of their stones $x$
that is adjacent to an opponent stone $y$ 
and {\em clobbers} $y$ with $x$
(removes $y$ and puts $x$ where $y$ was). 
Whoever cannot move loses. 

{\em Linear clobber} is clobber played on a path,
for example one row of a go board.
Figure~\ref{fdemo6} shows a linear clobber game from start (top row)
to finish (bottom row).
Each dot shows the stone that just clobbered.
Right/White played first, Left/Black won.

\begin{figure}[ht]
\hfill\includegraphics[scale=.5]{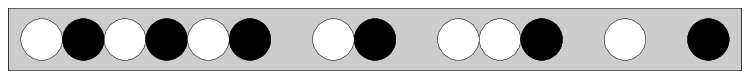}\hfill\ 
\caption{The parts from Figure~\ref{fdemo6}. The last two are trivial.}
\label{fparts}
\end{figure}

A {\em part} is a connected component of a clobber position.
For example, the position {\tt ox-oox} has two parts,
\texttt{ox} and \texttt{oox}.
$\mA$ is the set of non-empty even-length alternating-color parts
{\tt \{ox, oxox, oxoxox, ...\}}.
A part is {\em trivial} if it allows no legal move,
namely if each two adjacent stones are the same color.
$\mU$ is the set of all non-trivial parts that
appear in some game that starts from a part in $\mA$.
For example, each non-trivial part in Figure~\ref{fparts} is in $\mU$,
since each appears in the game in Figure~\ref{fdemo6}, 
which started from {\tt oxoxox} in $\mA$.
{\em Alternating Linear Clobber} (ALC) is 
the game of clobber starting from a game with all parts in $\mU$.

In clobber or more generally
any normal-play combinatorial game, 
whoever cannot move loses: there are no draws. 
Thus, each position is either
\begin{itemize}
\item {\em L-win} (L has a winning strategy when playing first) or 
\item {\em L-lose} (L has no winning strategy when playing second) and also
\item {\em R-win}  (R has a winning strategy when playing first) or 
\item {\em R-lose}  (R has no winning strategy when playing second).
\end{itemize}

Thus each position
is in exactly one outcome class:
$\mP$ ($L$-lose, $R$-lose),
$\mN$ ($L$-win,  $R$-win),
$\mL$ ($L$-win,  $R$-lose) or
$\mR$ ($R$-win,  $L$-lose).

In 2004, Albert et al.\ conjectured
that every part in $\mA$ except \verb+oxoxox+ is in $\mN$ \cite{AGNW}.
In this paper we prove their conjecture.

Here is how we found our proof.
In winter 2022, 
in his combinatorial games grad course,
Martin M\"{u}ller assigned students
the task of implementing a linear clobber solver:  
authors Folkersen and Randall
wrote solvers strong enough to solve 60-stone games.
(An enhanced version of this solver
now capable of solving 100-stone games
is now publicly available \cite{FBTM},
as is a more general solver \cite{FDM}.)
That summer, using solvers to generate data,
we looked for patterns that might lead
to a proof of the ALC conjecture.
From start game \texttt{a2n}, Left can move
to \texttt{o(2n-3)}:
this is left-right symmetric,
with only $n-2$ Right replies to consider.
In fall 2023 we pushed past this starting observation and
found the spiralling strategy described in this paper.

\section*{Combinatorial games background}
We use the following combinatorial games definitions and properties.
For proofs of these properties, see any relevant textbook such as
\cite{ANW16}, \cite{Haff} or \cite{CGTsiegel}.

A {\em game} $g$, written $g=\{G^L \: \vert \: G^R\}$,
is defined recursively by the set $G^L$ of {\em left options} (games that Left can move to) and the set
$G^R$ of {\em right options} (games that Right can move to).
For example, from {\tt oxoo} Left can move to 
{\tt x-oo} or {\tt o-xo} and Right can move to {\tt -ooo} or {\tt oo-o}
so {\tt oxoo} $=$ \{\{{\tt x-oo}, {\tt o-xo} \} \: $\vert$ \: 
    \{ {\tt -ooo}, {\tt oo-o}\}\}.
{\em Zero}, written 0, 
is the game $\{ \: \vert \: \}$ with no move options.

The {\em sum} of games 
$g=\{G^L \: \vert \: G^R\}$ and
$h=\{H^L \: \vert \: H^R\}$, written $g \myadd h$,
is defined recursively as the game
$j=\{J^L \: \vert \: J^R\}$,
where $J^L$ is the union of 
games $g_L \myadd h$ for all $g_L$ in $G^L$ and
games $g \myadd h_L$ for all $h_L$ in $H^L$.
$J^R$ is defined similarly.

The {\em negative}
of the game $g=\{G^L \: \vert \: G^R\}$
is the game $-g=\{-G^R \: \vert \:-G^L\}$,
where for a set of games $X$, $-X$ is the set $\{ -x \}$ for all $x$ in $X$.
To negate a clobber game, change the color of each stone.
For example, $-$({\tt xoxx}) is {\tt oxoo}.

Games $g$ and $h$ are {\em equivalent} if, for all games $x$,
$g \myadd x$ and $h \myadd x$ are in the same outcome class.
The sum of a game $g$ and a game $p$ in $\mP$ is equivalent to
the game $g$: within a sum, $p$ behaves like the game 0.
The sum of a game and its negative is in $\mP$: $p\myadd -p=0$.
We use this characterization:
$g$ and $h$ are equivalent if and only if $g \myadd -h$ is in $\mP$.

The {\em game tree} of a game is defined recursively
as the tree whose root is the game
and whose left subtrees are its Left options
and whose right subtrees are its Right options.
For example, the game tree of {\tt ox} has three nodes:
the root, its left child the game 0, and its right child the game 0.
The {\em canonical form} of a game $g$ is a simplest (fewest nodes
in its game tree) game $g'$ equivalent to $g$.
The canonical form of a game is unique.

Here is an exercise that uses this background material:
show that {\tt oxoo} is equivalent to the sum of {\tt xxo} and {\tt xo}
and that these games have canonical form \{0, {\tt xo} $\vert$ 0\}.

\section*{10 sets of parts}
For linear clobber,
define $\mA$ as the set of even-length alternating-color parts and $\mO$ as
the set of odd-length alternating-color parts that start and end with $o$. 
Building on $\mA$ and $\mO$, we define a total
of 10 sets of parts, shown below.
See Figure~\ref{fpartsten}.
Following convention, 
we use concatenation to describe paths,
e.g.\ {\tt o(ox)}$^3$ is $p=$ {\tt ooxoxox}.
We sometimes write a part as its 
initial monochromatic subtring,
followed by the number of stones in the part, 
optionally followed by its final 
monochromatic substring (if the final subtring has length at least two),
e.g.\ $p$ is oo7 and \texttt{ooxoxoo} is oo7oo.

{\centering

$\mA$ $=$ \verb+{ox, oxox, oxoxox, ...}+ $=$ \verb+{a2, a4, a6, ...}+ \\
$\mO$ $=$ \verb+{o, oxo, oxoxo, ...}+ $=$ \verb+{o, o3, o5, ...}+\\
$\moA$ $=$ \verb+{o, oox, ooxox,  ... }+ $=$ \verb+{o, oo3, oo5, ...}+\\
$\moO$ $=$ \verb+{oo, ooxo, ooxoxo, ...}+ $=$ \verb+{oo, oo4, oo6, ...}+\\
$\moOo$ $=$ \verb+{ooo, ooxoo, ooxoxoo, ...}+ $=$ \verb+{ooo, oo5oo, oo7oo, ...}+\\
$\mZ$ $=$ \verb+{ooxx, ooxoxx, ooxoxoxx, ...}+ $=$ \verb+{ooxx, oo6xx, oo8xx, ... }+\\
$\mX=-\mO=$   \verb+{x, xox, xoxox, ...}+ $=$ \verb+{x, x3, x5, ...}+\\
$\mAx=-\moA=$ \verb+{x, xxo, xxoxo, ... }+ $=$ \verb+{x, xx3, xx5, ...}+\\
$\mxX=-\moO=$ \verb+{xx, xxox, xxoxox, ...}+ $=$ \verb+{xx, xx4, xx6, ...}+\\
$\mxXx=-\moOo=$ \verb+{xxx, xxoxx, xxoxoxx, ...}+ $=$ \verb+{xxx, xx5xx, xx7xx, ...}+.

}

\begin{figure}[ht]
\hfill\includegraphics[scale=.5]{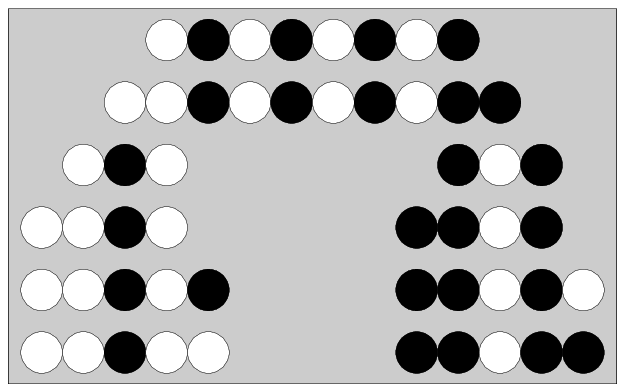}\hfill\ 
\caption{10 sets of parts.  From top: 
a8, oo10xx, (now left column) o3,
oo4, oo5, oo5oo,
(now right column) x3,
xx4, xx5, xx5xx.}
\label{fpartsten}
\end{figure}

\begin{theorem}
$\mU$ is the union of  
$\mA$, 
$\mO$, $\moA$, $\moO$, $\moOo$, $\mZ$,
$-\mO$, $-\moA$, $-\moO$ and $-\moOo$.
\label{thm0}
\end{theorem}

\begin{proof}[Proof of Theorem~\ref{thm0}.] 
$U=$ $\mA$ $\cup$ $\mO$ $\cup$ \ldots\ $\cup$ $\mZ$
is the set of paths that
can be written $\pm${\tt[o]}$a${\tt[c]}:
the path or its negative is 
an optional initial ${\tt o}$
followed by a path $a=$ {\tt (ox\ldots)} in $\mA$ or $\mO$
optionally followed by an extra {\tt c} ({\tt x} or {\tt o})
that is the same color as the last stone in $a$.
For example, Figure \ref{fallomoves} shows parts that can
arise when Right moves on a part in $U$.

To prove the theorem, 
it suffices to verify these properties:
\begin{itemize}
\item for every part $v$ resulting from a move on a part $u$ in $U$,
$v$ is in $U$,
\item for every part $u$ in $U$, some part $a$ in $\mA$
yields $u$ with at most two moves.
\end{itemize}

Consider a move on a part $p$ $=$ {\tt[o]}$a${\tt[c]}
in which stone $j$ clobbers neighboring stone $k$.
Since $j$ and $k$ have different colors, 
both $j$ and $k$ are in $a$.
This move corresponds to deleting $j$ and
color-switching $k$, making $k$ the end of its new path.
Thus, if $x$ and $y$ are same-color adjacent stones
in a part in $\mU$, then $x$ and/or $y$ 
is the end of its path.
\end{proof}

\begin{figure}[h]{\centering
{\small
\begin{tabular}{|c|c|c|c|} \hline 
$u$ & move & $a$ & $b$ \\ \hline
o($2t+1$) $\subseteq$ $\mO$ 
  & $2j+1$  $\rightarrow$ $2j+2$ & a($2j$) $\subseteq$  $\mA$ & oo($2t-2j$) $\subseteq$ $\moO$\\ \hline
x($2t+1$) $\subseteq$ $\mX$ 
  & $2j$  $\rightarrow$ $2j+1$ & x($2j-1$) $\subseteq$  $\mX$ & oo($2t-2j+1$) $\subseteq$ $\moA$\\ \hline
oo($2t+3$)o $\subseteq$ $\moOo$ 
  & $2j$  $\rightarrow$ $2j+1$ & oo($2j-1$) $\subseteq$  $\moA$ & oo($2t-2j+3$)oo $\subseteq$ $\moOo$\\ \hline
xx($2t+3$)xx $\subseteq$ $\mxXx$
  & $2j+1$  $\rightarrow$ $2j+2$ & xx($2j$) $\subseteq$  $Xx$ & a($2t-2j+2$) $\subseteq$ $\mA$\\ \hline
a($2t$) $\subseteq$ $\mA$ 
  & $2j+1$  $\rightarrow$ $2j$ & oo($2j$) $\subseteq$  $\moO$ & x($2t-2j-1$) $\subseteq$ $\mX$\\\hline
  & $2j+1$  $\rightarrow$ $2j+2$ & a($2j$) $\subseteq$  $\mA$ & oo($2t-2j-1$) $\subseteq$ $\moA$\\ \hline
oo($2t+2$)xx $\subseteq$ $\mZ$ 
  & $2j$  $\rightarrow$ $2j-1$ & oo($2j-1$)oo $\subseteq$  $\moOo$ & xx($2t-2j+2$) $\subseteq$ $\mxX$\\ \hline
  & $2j$  $\rightarrow$ $2j+1$ & oo($2j-1$) $\subseteq$  $\moA$ & oo($2t-2j+2$)xx $\subseteq$ $\mZ$\\ \hline
oo($2t+1$) $\subseteq$ $\moA$
  & $2j$  $\rightarrow$ $2j-1$ & oo($2j-1$)oo $\subseteq$  $\moOo$ & x($2t-2j+1$) $\subseteq$ $\mX$\\ \hline
  & $2j$  $\rightarrow$ $2j+1$ & xx($2j-1$) $\subseteq$  $\moA$ & xx($2t-2j+1$) $\subseteq$ $\moA$\\ \hline
oo($2t$) $\subseteq$ $\moO$
  &  $2j$  $\rightarrow$ $2j-1$ & oo($2j-1$)oo $\subseteq$  $oOo$ & a($2t-2j$) $\subseteq$ $A$\\ \hline
  &  $2j$  $\rightarrow$ $2j+1$ & oo($2j-1$) $\subseteq$  $oA$ & oo($2t-2j$) $\subseteq$ $oO$\\ \hline
xx($2t+1$) $\subseteq$ $\mAx$
  & $2j+1$  $\rightarrow$ $2j$ & oo($2j$)xx $\subseteq$  $\mZ$ & a($2t-2j$) $\subseteq$ $\mA$\\ \hline
  & $2j+1$  $\rightarrow$ $2j+2$ & xx($2j$) $\subseteq$  $\mxX$ & oo($2t-2j$) $\subseteq$ $\moO$\\ \hline
xx($2t+2$) $\subseteq$ $\mxX$ 
  & $2j+1$  $\rightarrow$ $2j$ & oo($2j$)xx  $\subseteq$  $\mZ$ & x($2t-2j+1$) $\subseteq$ $\mX$\\ \hline
  & $2j+1$  $\rightarrow$ $2j+2$ & xx($2j$) $\subseteq$  $\mxX$ & oo($2t-2j+1$) $\subseteq$ $\moA$\\ \hline
\end{tabular}
\caption{Right moves on part $u$ in $\mU$, yielding parts $a$ and $b$.}
\label{fallomoves}
}}\end{figure}

\section*{Our ALC Standard Form}
\begin{figure}[b]
\hfill\includegraphics[scale=.5]{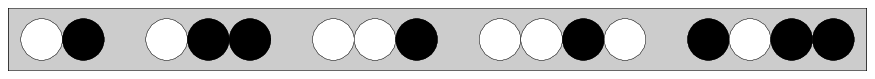}\hfill\ 
\caption{From left, games equivalent to $*$, \myup, \mydn, \myup$*$, \mydn$*$.}
\label{f5parts}
\end{figure}

The games in Figure~\ref{f5parts} are equivalent to the well known combinatorial games
star ($*$), up (\myup), down (\mydn), up-star (\myup$*$) and 
down-star (\mydn$*$). 
Zero (0) is any game with no legal moves.

\begin{figure}[ht]
\hfill\includegraphics[scale=.5]{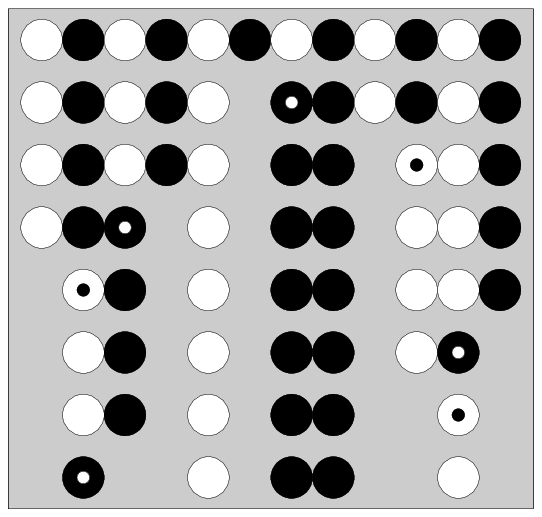}\hfill\ 
\caption{A game starting from {\tt ox12}. Left (black) plays first and wins.}
\label{fdemo12}
\end{figure}

In order to simplify proofs, we 
sometimes replace a game with an equivalent game.
For example, consider {\tt (ox)}$^6$ $=$ {\tt a12} $=$ {\tt oxoxoxoxoxox}: 
see Figure~\ref{fdemo12}.
In our proofs, we replace {\tt (ox)}$^6$ 
with the equivalent game {\tt oxox}$*$.
Algorithm ASF converts an arbitrary ALC game into
our standard form. Its correctness relies on these properties:
\begin{itemize}
\item a game's outcome is unchanged if a part 
is replaced with an equivalent part,
\item a game with no legal moves is game 0,
\item the sum of a game and its negative is equivalent to the game 0.
\end{itemize}

The {\em negative} of each of the following rules is obtained
by replacing each part of the left-hand side and right-hand side
with its negative. 
For example, the negative of rule $\alpha$ 
is to replace each of {\tt x, xx, xxx, xxoo, xoxoxo} with the game 0.

{\bf Algorithm Standard Form (ASF)} 
{\em
Given any ALC game, repeatedly apply each rule in preference order 
$(\alpha,-\alpha,\beta,-\beta,\ldots,\eta,-\eta)$
until no rule applies.
}

{\tt
\begin{tabular}{rl}
$\alpha$. & o, oo, ooo, ooxx, oxoxox $\Rightarrow$ \verb+-+ (game with no stones) \\
$\beta$. & for any parts $p$ and $-p$, $p \myadd -p \Rightarrow$ \verb+-+ \\
$\gamma$. & oxo, ooxox, ooxoxoo, xxoxoxx $\Rightarrow$ ox \\
$\delta$. & oxoxoxoxo $\Rightarrow$ ooxo \\
$\varepsilon$. & ooxoxx, oxoxoxoxoxox $\Rightarrow$ oxox $\myadd$ ox \\
$\zeta$. & ooxoo $\Rightarrow$ oox \\
$\eta$. & ooxo $\Rightarrow$ xxo $\myadd$ ox
\end{tabular}
}

\begin{theorem}
ASF terminates and returns a game equivalent to the input game.
\label{thmcorrect}
\end{theorem}

\begin{proof}[Proof of Theorem~\ref{thmcorrect}.]
Games $g$ and $h$ are equivalent if and only if $g \myadd -h$ is in $\mP$.
The reader can verify 
by hand or with a software package such as CGSuite \cite{cgsuite}
that each rule yields an equivalent game:
proving by hand the correctness of Rule $\varepsilon$ requires some effort.
Each rule reduces the sum, over all parts, of the square of the number of stones
in that part: thus the algorithm terminates.
We need this sum-of-squares condition only because of Rule $\eta$: 
every other rule reduces the total number of stones.
\end{proof}

Here are sets $\mA$ to $\mZ$ again, this time in standard form \ldots

\begin{tabular}{cccccccccccc}
$\mA$ =   & \{ &   & {\tt a2} &{\tt a4} & 0 & {\tt a8} & {\tt a10}
           & {\tt a4 $\myadd$ a2} & {\tt a14} & \ldots\}\\
$\mO$ =   & \{ & 0 & {\tt a2} & {\tt o5} & {\tt o7} & {\tt o9} & {\tt o11} &  {\tt o13} & {\tt o15}  &\ldots\}\\
$\moO$ =  & \{ && 0 & {\tt xxo $\myadd$ a2} & {\tt oo6} & {\tt oo8} & {\tt oo10} & {\tt oo12}  & {\tt oo14} &\ldots\}\\
$\moA$ =  & \{ & 0 & \texttt{oo3} & {\tt a2} & {\tt oo7} & {\tt oo9} & {\tt oo11} & {\tt oo13} & {\tt oo15} &\ldots\}\\
$\moOo$ = & \{ && 0 & \texttt{oo3} & {\tt a2} & {\tt oo9oo} & {\tt oo11oo} & {\tt oo13oo} & {\tt oo15oo}  &\ldots\}\\
$\mZ$ = & \{ &&& 0 & {\tt a4 $\myadd$ a2} & {\tt oo8xx} & {\tt oo10xx} & {\tt oo12xx} & {\tt oo14xx} &\ldots\}\\
\end{tabular}

\ldots\ and here are five respective subsets $\mA'$ to $\moOo'$, 
which are central to our main results:

\begin{tabular}{ccccccccccccl}
$\mA'$ =  & \{ & \phantom{0} & \phantom{0} & &  & {\tt a8} & {\tt a10}
           & & {\tt a14} & \ldots\}\\
$\mO'$ =   & \{ &   &     & {\tt o5} & {\tt o7} & & {\tt o11} &  {\tt o13} & {\tt o15}  &\ldots\}\\
$\moO'$ =  & \{ & & & & & & & {\tt oo10} & {\tt oo12}  &\ldots\}\\
$\moA'$ =  & \{ & & & & {\tt oo7} & {\tt oo9} & {\tt oo11} & {\tt oo13} & {\tt oo15} &\ldots\}\\
$\moOo'$ = & \{ & & & \texttt{oo3} & & {\tt oo9oo} & {\tt oo11oo} & {\tt oo13oo} & {\tt oo15oo}  &\ldots\} & .\\
\end{tabular}

\section*{Main results}
We need a last bit of terminology.
Let

{\centering
$\mA' = \{{\tt a8}, {\tt a10}, {\tt a14}, \ldots\}$, ~  
$\mO' = \{{\tt o5}, {\tt  o7}, {\tt o11}, \ldots\}$, ~ 
$\moO' = \{{\tt oo10}, {\tt  oo12}, \ldots\}$,

$\moA' = \{{\tt oo7}, {\tt  oo9}, \ldots\}$, ~
$\moOo' = \{{\tt oo3}, {\tt  oo9oo}, {\tt oo11oo}, \ldots\}$, ~
$\mI$ $=$ $\{{\tt a2}, {\tt a4}, {\tt oo6}\}$,

$Q$ $=$ \{{\tt o5} $\myadd$ {\tt a2},
~ {\tt o7} $\myadd$ {\tt a4}, 
~ {\tt o15} $\myadd$ {\tt a4} $\myadd$ {\tt a2}\}.

}

$\mI$ is a set of irregular games that we handle separately in our proofs.
$Q$ is a set of games that Left needs to avoid leaving for Right.
The other sets above are defined in the previous section.
For any game $g$ in standard form 
and with all parts in 

{\centering
$K=\mO' \cup \moO' \cup \moOo' \cup \mI \cup \{{\tt xxo}\} \cup
\{{\tt oo8}\} \cup \mA' \cup \moA'$

}

we define the {\em count vector} $\mu(g)=$
$(a,b,c,d,e,f,y,z)$, where $a \ldots y,z$ are the numbers 
of not necessarily distinct parts of $g$ that are
in the respective seven sets of the above union.
When $g$ has no parts in $\mA'$ nor $\moA'$
we write this vector as $(a,b,c,d,e,f)$ instead of
$(a,b,c,d,e,f,0,0)$. See Figure~\ref{fSeg}.
Let $\mS$ be the set of games $g$ in $K$ 
with at least one part and
with count vector $(a,\ldots, z)$ such that

\begin{enumerate}\addtocounter{enumi}{-1}
\item ($y,z$) is (1,0) and $a \geq c$, or 
($y,z$) is (0,0) and $a \geq c$, or 
($y,z$) is (0,1) and $a \geq c+1$,
\item $(y,z)$ is (0,0) and $a\geq c+1$ and 
the game $g$ is not in $Q$,
\item $(y,z,a,b,c)$ is (0,0,0,0,0)
and $e\geq 1$ and either $d=0$ or $d+e\geq 3$.
\end{enumerate}
$\mS_0, \mS_1, \mS_2$ 
are the subsets of $\mS$ that satisfy properties 0,1,2 respectively. 
Notice that $\mS_1$ and $\mS_2$ are non-intersecting subsets of $\mS_0$.
$LL$ is the set
\{oo8, ~ o5 $\myadd$ oox $\myadd$ a4 $\myadd$ a2, \\
a14 $\myadd$ xxo $\myadd$ a4 $\myadd$ a2\}.
Each game in $LL$ is in $\mL$. 

\begin{figure}[h]{\centering
{\small
\begin{tabular}{|c|c|c|} \hline
  game & count vector & subset of $\mS$ \\ \hline
{\tt o11} $\myadd$
{\tt oo10} $\myadd$
{\tt oo9oo} $\myadd$
{\tt a2} $\myadd$
{\tt a8}
&
(1,1,1,1,0,0,1,0) & $\mS_0$ \\ \hline
{\tt o11} $\myadd$
{\tt oo14} $\myadd$
{\tt oo9oo} $\myadd$
{\tt a2} $\myadd$
{\tt oo9}
&
(1,1,1,0,1,0,0,1) & $\mS_0$ \\ \hline
{\tt oo8} $\myadd$ {\tt oo14}&
(0,1,0,0,0,1,0,0) & $\mS_0$ \\ \hline
{\tt o5} $\myadd$
{\tt o11} $\myadd$
{\tt oo6} $\myadd$
{\tt oo9oo} $\myadd$
{\tt a2}&
(2,1,1,1,0,0,0) & $\mS_1$ \\ \hline
{\tt oo6} $\myadd$
{\tt a4} $\myadd$
{\tt xxo} $\myadd$
{\tt xxo} 
& (0,0,0,2,2,0) & $\mS_2$ \\ \hline
{\tt a4} $\myadd$ 
{\tt a2} $\myadd$
{\tt xxo} & 
(0,0,0,2,1,0) & $\mS_2$ \\ \hline
{\tt xxo} & (0,0,0,0,1,0) & $\mS_2$ \\ \hline
\end{tabular}

}}
\caption{Example games in $\mS$.}
\label{fSeg}
\end{figure}

Here are our results.

\begin{theorem}
On any game in $\mS_1\cup\mS_2$, 
Right can move only to $\mS_0$.
\label{thmRtoS0}
\end{theorem}

\begin{theorem}
On any game in $\mS_0$,
Left can move to a game in $\mS_1 \cup \mS_2 \cup LL \cup \{0\}$.
\label{thmLonS0}
\end{theorem}

\begin{corollary*}
$\mS_0$ is in $\mL \cup \mN$.
$\mS_1 \cup \mS_2$ is in $\mL$.
Each game in $\mA$ $\setminus$ \{{\tt oxoxox}\} is in $\mN$.
\end{corollary*}

\begin{figure}[h]
\hfill\includegraphics[scale=.6]{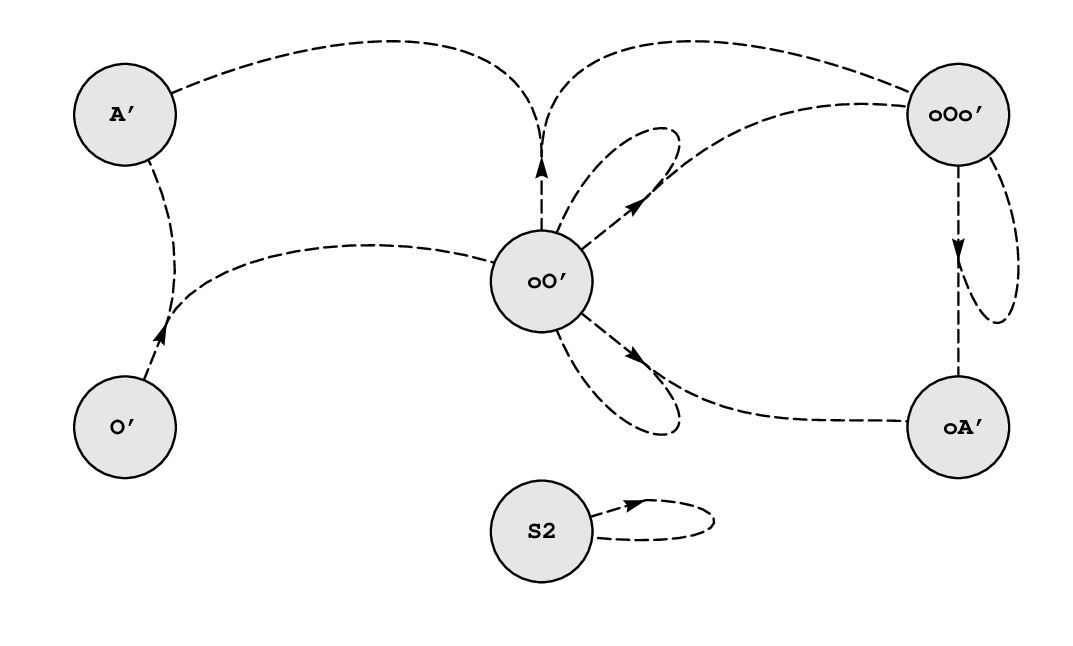}\hfill\ 
\caption{How Right can move on $\mS$.}
\label{fdiagRight}
\end{figure}

\begin{figure}[h]
\hfill\includegraphics[scale=.6]{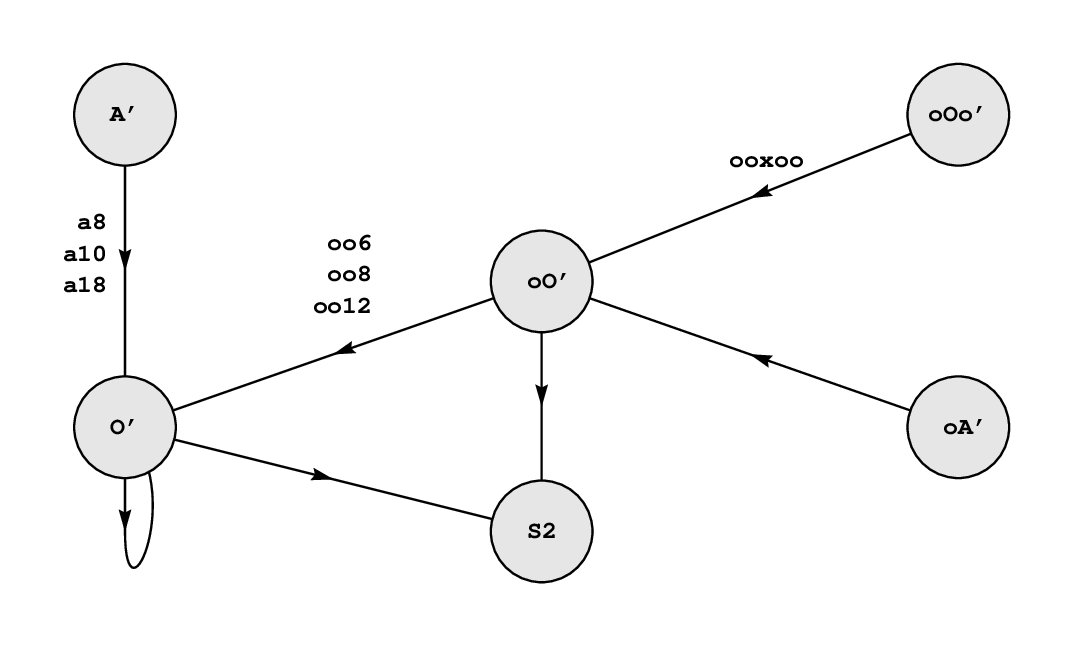}\hfill\ 
\caption{How Left plays on $\mS$. Edge labels are exceptions.}
\label{fdiagLeft}
\end{figure}

\begin{figure}[h]
\hfill\includegraphics[scale=.6]{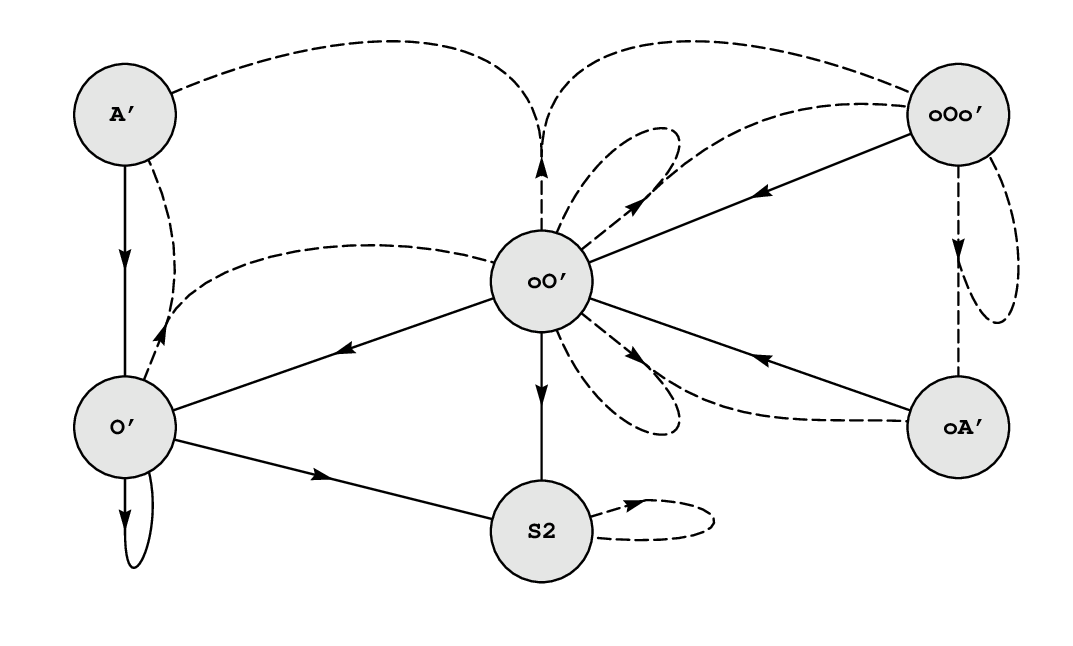}\hfill\ 
\caption{How Left (solid arrows) and Right (dashed) play on $\mS$.}
\label{fdiagram}
\end{figure}

The last statement of the corollary is the conjecture of Albert et al.
Figures~\ref{fdiagRight} through \ref{fdiagram}
show how Right can play on $\mS$ and how Left plays
on $\mS_0$ in the proof of Theorem~\ref{thmLonS0}.

\section*{Proofs}
We begin with a motivating observation.

\begin{obs*}
The sum of a game in $\mL$ and 
a game in $\mL \cup \mN$ is in $\mL \cup \mN$.
\end{obs*}

\begin{proof}[Proof of Observation.]
Let $h$ and $k$ be in $\mL$ and $\mL \cup \mN$ respectively.
We want to show that $h \myadd k$ is $L$-win so assume $L$ plays first.
$L$ wins by following simultaneously
a second-player-win strategy on $h$
and a first-player-win strategy on $k$.
\end{proof}

\begin{proof}[Proof of Theorem \ref{thmRtoS0}.]
Assume that Right plays
on part $p$ of a game $g$ in $\mS_1 \cup \mS_2$, 
leaving game $h=g-p \myadd j \myadd k$,
where $j$ and $k$ are the parts that arise from playing
on $p$: one or both of $j,k$ might be 0. 
We then apply ASF to $h$, yielding $h$'.
We want to show that $h'$ is in $S_0$.
Let $(a,b,\ldots,y,z)$ and
$(a',b',\ldots,y',z')$ 
be the count vectors of $g$ and $h'$ respectively.

Case 1: $g$ is in $S_1$.
~ Thus $g$ has 
~ $a\geq c+1$ parts from $\mO'$,
~ $b$ parts from $\moO'$,
~ $c$ parts from $\moOo'$,
~ $d$ parts from $\mI$,
~ $e$ parts from \{\texttt{xxo}\},
~ $f$ parts from \{\texttt{oo8}\}, and 
~ no parts from $\mA'$ nor $\moA'$.
~ There are 6 subcases, depending on whether
$p$ is in $\mO'$, $\moO'$, \ldots, or \{\texttt{oo8}\} respectively.
Figures \ref{fRmvO'}, \ref{fRmv2}, and \ref{fRmvoOo'}
show the changes to the count vector that arise
after a move on a part in $\mO'$, $\moO'$, and $\moOo'$ respectively.

\begin{figure}{\centering
{\small \begin{tabular}{|c|c|c|c|c|c|c|c|c|c|} \hline 
\multicolumn{10}{|c|}{Right move on $\mO'$ yields parts $j,k$ and change to vector}\\ \hline
$j$ & $k$ & $\Delta_a$ & $\Delta_b$ & $\Delta_c$ & $\Delta_d$ & $\Delta_e$ & $\Delta_f$& $\Delta_y$ & $\Delta_z$ \\ \hline
$\emptyset$ & \texttt{a2} & -1 & 0 & 0 & +1 & 0 & 0 & 0 & 0 \\ \hline
$\emptyset$ & \texttt{a4} & -1 & 0 & 0 & +1 & 0 & 0 & 0 & 0 \\ \hline
$\emptyset$ & $\mA'$ & -1 & 0 & 0 & 0 & 0 & 0 & +1 & 0 \\ \hline
\texttt{oo4} = \texttt{xxo} + a2 & $\emptyset$ & -1 & 0 & 0 & +1 & +1 & 0 & 0 & 0 \\ \hline
\texttt{oo4} = \texttt{xxo} + a2 & \texttt{a2} & -1 & 0 & 0 & +2 & +1 & 0 & 0 & 0 \\ \hline
\texttt{oo4} = \texttt{xxo} + a2 & \texttt{a4} & -1 & 0 & 0 & +2 & +1 & 0 & 0 & 0 \\ \hline
\texttt{oo4} = \texttt{xxo} + a2 & $\mA'$ & -1 & 0 & 0 & +1 & +1 & 0 & +1 & 0 \\ \hline
\texttt{oo6} & $\emptyset$ & -1 & 0 & 0 & +1 & 0 & 0 & 0 & 0 \\ \hline
\texttt{oo6} & \texttt{a2} & -1 & 0 & 0 & +2 & 0 & 0 & 0 & 0 \\ \hline
\texttt{oo6} & \texttt{a4} & -1 & 0 & 0 & +2 & 0 & 0 & 0 & 0 \\ \hline
\texttt{oo6} & $\mA'$ & -1 & 0 & 0 & +1 & 0 & 0 & +1 & 0 \\ \hline
\texttt{oo8} & $\emptyset$ & -1 & 0 & 0 & 0 & 0 & +1 & 0 & 0 \\ \hline
\texttt{oo8} & \texttt{a2} & -1 & 0 & 0 & +1 & 0 & +1 & 0 & 0 \\ \hline
\texttt{oo8} & \texttt{a4} & -1 & 0 & 0 & +1 & 0 & +1 & 0 & 0 \\ \hline
\texttt{oo8} & $\mA'$ & -1 & 0 & 0 & 0 & 0 & +1 & +1 & 0 \\ \hline
$\moO'$ & $\emptyset$ & -1 & +1 & 0 & 0 & 0 & 0 & 0 & 0 \\ \hline
$\moO'$ & \texttt{a2} & -1 & +1 & 0 & +1 & 0 & 0 & 0 & 0 \\ \hline
$\moO'$ & \texttt{a4} & -1 & +1 & 0 & +1 & 0 & 0 & 0 & 0 \\ \hline
$\moO'$ & $\mA'$ & -1 & +1 & 0 & 0 & 0 & 0 & +1 & 0 \\ \hline
\end{tabular}}

}\caption{For game $g$, $R$ moves on part $p$ in $\mO'$.}
\label{fRmvO'}
\end{figure}

Subcase 1.1: $p$ is in $\mO'$ $=$ \{\texttt{o5, o7, o11,} \ldots\}, so
we can assume that $j,k$ are in sets
\{0, \texttt{a2}, \texttt{a4}, \ldots \} and
\{0, \texttt{oo4}, \texttt{oo6}, \ldots \}, so sets
\{0, \texttt{a2, a4}\} $\cup$ $\mA'$ and 
\{0, \texttt{xxo+a2, oo6, oo8}\} $\cup$ $\moO'$ respectively.

Recall that $y,z$ are both 0.  We have
either $y' = y = 0$ or  $y'= y+1 = 1$ (whenever $j$ is in $\mA'$),
and $z' = z = 0$, and
$a' \geq c'$ (because $h'$ does not include $p$ which is in $\mO$ 
and neither $j$ nor $k$ is in $\mO'$ nor $\moOo'$, 
so $a' = a - 1$ and $c'=c$).
It follows that $h'$ is in $S_0$ $\cup$ \{0\}.
To finish the proof in Subcase 1.1 it remains to show that $h'$ is not 0, namely
$h'$ has at least one part.
This holds if $k$ is \texttt{oo2t} where $2t \geq 6$, so assume
$k$ is \texttt{oo} or \texttt{ooxo}.

In the former case $k$ is 0 so we have $p$ in $g$ replaced
by $j$ in $\mA$ and $h'$ has a part in $\mA'$ unless
$p$ is \texttt{o5}, \texttt{o7}, or \texttt{o15}
in which $j$ is \texttt{a2}, \texttt{a4}, or
\texttt{a12} is equivalent to \texttt{a4 + a2} respectively.
But then the only way that these parts
are not in $h'$ is when $g$ is \texttt{o5+a2},
\texttt{o7+a4}, or \texttt{o15+a4+a2} respectively,
contradicting our assumption that $g$ is not in $Q$.

In the latter case $k$ is \texttt{ooxo} and
so will be replaced with the 
equivalent \texttt{xxo+ox} in $h'$. But for \texttt{xxo} 
to then be eliminated by negation from $h'$,
it must be that \texttt{oox} is a part of $g$, 
so $c \geq 1$, so $a\geq 2$
and since R plays on $p$ in $\mO'$ we have
$a' = a -1 \geq 2-1=1$ so $h'$ has at least one part and we are done.

Subcase 1.2: $p$ is in $\moO'$ $=$ \{\texttt{oo10, oo12, \ldots}\}.
Here $j,k$ are $\moA,\moO$ respectively or $\moOo,\mA$ respectively.
Consider the former case: $j,k$ are in sets
\{0, \texttt{oox}, \texttt{ooxox}, \ldots \} and
\{0, \texttt{ooxo}, \texttt{ooxoxo}, \ldots \}, so sets
\{0, \texttt{oo3, a2}\} $\cup$ $\moA'$ and 
\{0, \texttt{xxo+a2, oo6, oo8}\} $\cup$ $\moO'$ respectively.
Recall that $y,z$ are both 0.  We have
$y' = y = 0$, and
either $z' = z = 0$ or $z' = z + 1 = 1$ (whenever $j$ is in $\moA'$),
When $z'=0$, $c'\leq c+1$ so $a' \geq c'$.
When $z'=1$, $c'\leq c$ so $a' \geq c'+1$ and we are done in the former case.

Next consider the latter case: $j,k$ are in sets
\{0, \texttt{ooxoo}, \texttt{ooxoxoo}, \ldots \} and
\{0, \texttt{ox}, \texttt{oxox}, \ldots \}, so sets
\{0, \texttt{oo3, a2}\} $\cup$ $\moOo'$ and 
\{0, \texttt{a2, a4}\} $\cup$ $\mA'$ respectively.
Again, $y,z$ are both 0.  We have
either $y' = y = 0$ or $y' = y + 1 = 1$ (whenever $k$ is in $\mA'$),
and $z' = z = 0$, and
$a' \geq c'$ (because $h'$ does not include $p$ which is in $\mO$ 
and neither $j$ nor $k$ is in $\mO'$ and $j$ is not in $\moOo'$).
It follows that $h'$ is in $S_0$ $\cup$ \{0\} 
and we are done in the latter case.

To finish the proof in Subcase 1.2
it remains to show that $h'$ is not 0, namely $h'$ has at least one part.
Since $a\geq c+1$, $g$ 
has at least one part in $\mO'$, and this part is not eliminated 
by negation by a subgame arising from Right's move on $p$ in $\moOo'$. 
So $h'$ has at least one part and we are done.

Subcase 1.3: $p$ is in $\moOo'$ $=$ \{\texttt{oo3, oo9oo, oo11oo, \ldots}\}.
Here $j,k$ are in sets 
\{0, \texttt{oo5oo = oo3}, \texttt{oo7oo = a2}, \texttt{oo9oo}, \ldots \} and
\{0, \texttt{oo3}, \texttt{oo5 = a2}, \texttt{oo7}, \ldots \}, so sets
\{0, \texttt{oo3, a2}\} $\cup$ $\moOo'$ and \{0, \texttt{oo3, a2}\} $\cup$ $\moA'$ respectively.
Recall that $a \geq c+1$ and 
that $y,z$ are both 0.  We have
$y' = y = 0$ and
either $z' = z = 0$ or $z'= z+1 = 1$ (whenever $k$ is in $\moA'$).
When $z'=0$, $c'\leq c+1$ so $a' \geq c'$.
When $z'=1$, $c'\leq c$ so $a' \geq c'+1$:
$p$ is in $\moOo'$, $j$ might be in $\moOo'$, and $k$ is not in
$\moOo'$, so $c' \leq c+1$.
It follows that $h'$ is in $S_0$ $\cup$ \{0\}.

To finish the proof in Subcase 1.3 
it remains to show that $h'$ is not 0, namely $h'$ has at least one part.
Since $a\geq c+1$, $g$ 
has at least one part in $\mO'$, and this part is not negated by 
any subgame arising from Right's move on $p$ in $\moOo'$. 
So $h'$ has at least one part and we are done.

Subcase 1.4: $p$ is in $\mI$.
Each possible move by $R$ leaves all parts in 
\{0, \texttt{a2}, \texttt{oo3}, \texttt{xxo}, 
\texttt{xxo} $\myadd$ \texttt{a2} \}
and again we have $a_h\geq c_h$ and we are done.

Subcase 1.5: 
$p$ is in \{\texttt{xxo}\}.
Now R's only move is to
\{\texttt{xo}\} which is in $\mI$ so, like $g$, $h$ is in $S_0$.
Also, since $p$ is not in $\mO'$ nor in $\moOo'$ we have
$a=a_h$ and $c=c_h$, and since
$a \geq c+1$ we have $a_h\geq c_h+1$,
so $h$ is in $S_1$ and so $S_0$ and we are done.

Subcase 1.6: $p$ is in \{\texttt{oo8}\}.
Each possible move by $R$ leaves all parts in $\mI$
or
\texttt{ooxoo}$\myadd$\texttt{xo} $=$ \texttt{oox}$\myadd$\texttt{ox}.
In all cases we have $a_h\geq c_h+1 - 1$, so $a_h\geq c_h$,
so $h$ is in $S_0$ and we are done.
This concludes the proof in Case 1.
\vspace*{5pt}

\begin{figure}{\centering
{\small \begin{tabular}{|c|c|c|c|c|c|c|c|c|c|} \hline 
\multicolumn{10}{|c|}{Right move on $\moO'$ yields parts $j,k$ and change to vector}\\ \hline
$j$ & $k$ & $\Delta_a$ & $\Delta_b$ & $\Delta_c$ & $\Delta_d$ & $\Delta_e$ & $\Delta_f$& $\Delta_y$ & $\Delta_z$ \\ \hline
$\emptyset$ & $\mA'$ & 0 & -1 & 0 & 0 & 0 & 0 & +1 & 0 \\ \hline
$\emptyset$ & \texttt{a2} & 0 & -1 & 0 & +1 & 0 & 0 & 0 & 0 \\ \hline
$\emptyset$ & \texttt{a4} & 0 & -1 & 0 & +1 & 0 & 0 & 0 & 0 \\ \hline
$\emptyset$ & \texttt{oo4} = \texttt{xxo} + a2 
& 0 & -1 & 0 & +1 & +1 & 0 & 0 & 0 \\ \hline
$\emptyset$ & \texttt{oo6} & 0 & -1 & 0 & +1 & 0 & 0 & 0 & 0 \\ \hline
$\emptyset$ & \texttt{oo8} & 0 & -1 & 0 & 0 & 0 & +1 & 0 & 0 \\ \hline
$\emptyset$ & $\moO'$ & 0 & 0 & 0 & 0 & 0 & 0 & 0 & 0 \\ \hline
\texttt{oo3} $\subseteq \moOo'$ & $\emptyset$ & 0 & -1 & +1 & 0 & 0 & 0 & 0 & 0 \\ \hline
\texttt{oo3} $\subseteq \moOo'$ & \texttt{a2} & 0 & -1 & +1 & +1 & 0 & 0 & 0 & 0 \\ \hline
\texttt{oo3} $\subseteq \moOo'$ & \texttt{a4} & 0 & -1 & +1 & +1 & 0 & 0 & 0 & 0 \\ \hline
\texttt{oo3} $\subseteq \moOo'$ &  \texttt{oo4} = \texttt{xxo} + a2 
& 0 & -1 & -1 & +1 & +1 & 0 & 0 & +1 \\ \hline
\texttt{oo3} $\subseteq \moOo'$ & $\moO'$ 
& 0 & 0 & +1 & 0 & 0 & 0 & 0 & 0 \\ \hline
$\moA'$ & $\emptyset$ & 0 & -1 & 0 & 0 & 0 & 0 & 0 & +1 \\ \hline
$\moA'$ & \texttt{oo4} = \texttt{xxo} + a2 & 0 & -1 & 0 & +1 & +1 & 0 & 0 & +1 \\ \hline
$\moA'$ &  \texttt{oo6} & 0 & -1 & 0 & +1 & 0 & 0 & 0 & +1 \\ \hline
$\moA'$ & \texttt{oo8} & 0 & -1 & 0 & 0 & 0 & +1 & 0 & +1 \\ \hline
$\moA'$ & $\moO'$ & 0 & 0 & 0 & 0 & 0 & 0 & 0 & +1 \\ \hline
$\moOo'$ & $\emptyset$ & 0 & -1 & +1 & 0 & 0 & 0 & 0 & 0 \\ \hline
$\moOo'$ & $\mA'$ & 0 & -1 & +1 & 0 & 0 & 0 & +1 & 0 \\ \hline
\end{tabular}}

}\caption{For game $g$, Right moves on part $p$ in $\moO'$.}
\label{fRmv2}
\end{figure}

\begin{figure}{\centering
{\small \begin{tabular}{|c|c|c|c|c|c|c|c|c|c|} \hline 
\multicolumn{10}{|c|}{Right move on $\moOo'$ yields parts $j,k$ and change to vector}\\ \hline
$j$ & $k$ & $\Delta_a$ & $\Delta_b$ & $\Delta_c$ & $\Delta_d$ & $\Delta_e$ & $\Delta_f$& $\Delta_y$ & $\Delta_z$ \\ \hline
$\emptyset$ & $\moA'$ & 0 & 0 & -1 & 0 & 0 & 0 & 0 & +1 \\ \hline
$\moOo'$ & $\emptyset$ & 0 & 0 & 0 & 0 & 0 & 0 & 0 & 0 \\ \hline
$\moOo'$ & \texttt{a2} & 0 & 0 & 0 & +1 & 0 & 0 & 0 & 0 \\ \hline
$\moOo'$ & \texttt{oo3} $\subseteq \moOo'$ & 0 & 0 & +1 & 0 & 0 & 0 & 0 & 0 \\ \hline
$\moOo'$ & $\moA'$ & 0 & 0 & 0 & 0 & 0 & 0 & 0 & +1 \\ \hline
\end{tabular}}

}\caption{For game $g$, Right moves on part $p$ in $\moOo'$.}
\label{fRmvoOo'}
\end{figure}

\begin{figure}{\centering
{\small \begin{tabular}{|c|c|} \hline 
\multicolumn{2}{|c|}{Right move on $p$ in $S_2$ yields 
new part(s) $j$} \\ \hline
$p$ & $j$ \\ \hline
\texttt{a2} & 0 \\ \hline
\texttt{a4} & \texttt{a2 or oox} \\ \hline
\texttt{oo6} & \texttt{a2 or oox or a2+xxo} \\ \hline
\texttt{xxo} & \texttt{a2} \\ \hline
\texttt{oo8} & \texttt{a2, a4, oo6, or a2+oox} \\ \hline
\end{tabular}}

}\caption{For game $g$, Right moves on part $p$ in $S_2$.}
\label{fRmvS2}
\end{figure}

Case 2. $g$ is in $S_2$ with count vector $(0, 0, 0, d, e, f)$.
Subcase 2.1: $d=0$ and $e \geq 1$.
Subcase 2.2: $d \geq 1$ and $d+e \geq 3$. The following arguments apply to both subcases.

Each part of $g$ is in 
$S_3$ $=$ \{\texttt{a2, a4, oo6, xxo, oo8}\}.
Figure~\ref{fRmvS2} shows how Right can move on
each possible part $p$ of $g$,
yielding part(s) $j$ of the resulting
game $h$ that is then
replaced by the equivalent game $h'$ by ASF.
Because $g$ has $d\geq 1$, 
whenever $j$ is \texttt{oox} or
\texttt{a2+oox}, \texttt{oox}
negates a part \texttt{xxo} in $g$,
so $h'$ has no part \texttt{oox},
so each part in $h'$ is $S_3$,
so $h'$ is in $S_0 \cup \{0\}$. 
It remains only to show that $h'$ is not 0,
namely has at least one part.

As in Figure~\ref{fRmvS2},
consider each possible $p$ that Right moves on.
First consider $p=$ \texttt{a2}.
Right plays \texttt{a2} to 0, so $j$ is 0.
Here $h = g - p + j = g -$\texttt{a2}$+ 0$,
so after applying ASF to $h$ we have $h' = g-$\texttt{a2}.
Every game in $S_2$ has at least at least one part \texttt{xxo},
so here $g$ a part \texttt{xxo}, so $h'$ does as well and we are done.

Next consider $p=$ \texttt{a4}.
Thus $g$ includes \texttt{a4} and is in $S_2$,
so $g$ also includes \texttt{xxo} and has at least three parts
in $\mI'$ \{\texttt{xxo}\} $\cup$ $\mI$.
Right plays \texttt{a4} to $j$ $=$ \texttt{a2} or $j$ $=$ \texttt{oox},
which erases at most one part of $\mI'$,
so at least one part of $\mI'$ remains in $h'$ and we are done.
By similar arguments, for each $p$ in \{\texttt{oo6, xxo}\},
$h'$ has at least one part in $\mI'$ so we are done.

So finally assume that $p$ is \texttt{oo8}:
$g$ includes \texttt{oo8} and is in $S_2$,
so either $g$ has at least one part in \{\texttt{xxo}\}
(Subcase 2.1)
or at least three parts in \{\texttt{xxo}\} $\cup$ $\mI$
(Subcase 2.2).
By Figure~\ref{fRmvS2},
Right plays \texttt{oo8} to $j$ $=$ \texttt{a2}
or \texttt{a4}
or \texttt{oo6}
or \texttt{a2 + oox},
erasing from $\mI'$ at most one part in Subcase 2.1
(because here $\mI'=\mI$)
and at most two parts of in Subcase 2.2.
In the former case, if a part of $\mI'$ is erased
then it must be part \texttt{xxo} was negated by
adding parts \texttt{a2} and \texttt{oox},
so the part \texttt{a2} remains in $h'$.
In the later case,
$\mI'$ started with at least three parts and at most
two were negated, so again at least one such part
remains in $h'$ and we are done.
This concludes the proof of Case 2 and
the proof of Theorem~\ref{thmRtoS0}.
\end{proof}

\begin{proof}[Proof of Theorem \ref{thmLonS0}.]
Assume that Left plays on a game
$g$ with count vector $(a,b,c,d,e,f,y,z)$
in $\mS_0$ by following Rules 1 through 7 as shown
in Figures~\ref{fr1} through \ref{fr67},
leaving game $h$.
We want to show that $h$ is in 
$\mS_1 \cup \mS_2 \cup LL \cup \{0\}$.
Left considers Rules 1 through 7 in order.

\begin{figure}[h]{\centering
\input{r1}

}\caption{How Left plays on $\mS_0$: rule 1.}
\label{fr1}
\end{figure}

Subcase 1: Rule 1 applies. Thus $\mA'$ is non-empty,
so $y=1$ ($y$ is positive and $g$ is in $S_0$).
If any of Rule 1a,1b,1c applies then Left plays as shown
in Figure~\ref{fr1}: for example with Rule 1a play
\texttt{a8+a2} to \texttt{xxo+a4+a2}.
We leave it to the reader to check -- for example with CGSuite --
that each resulting game $h$ is in $\mL$,
in which case we are done.
If Rule 1d applies then Left plays on part 
\texttt{a8} to \texttt{o5}, or
\texttt{a10} to \texttt{o7}, and so on,
in which case the count vector $(a',b',c',d',e',f',y',z')$
of $h'$ satisfies $y'=0$ 
($y=1$ and Left played a part in $\mA'$ to a part in $\mO'$)
and $a'\geq c'+1$ 
(the number of parts in $\mO'$ increased,
the number of parts in $\moOo'$ did not change),
so $h'$ is in $\mS_1$ and we are done.

\begin{figure}[h]{\centering
\input{r2}

}\caption{How Left plays on $\mS_0$: rule 2.}
\label{fr2}
\end{figure}

Subcase 2: Rule 1 does not apply (so $y=0$),
Rule 2 applies.
Thus  $z=1$, since $g$ is in $S_0$ and $\moA'$ is non-empty.
Left plays as shown in Figure~\ref{fr2},
so $z'=z-1=0$ and $b'=b$ 
(when Left's move yields 
\texttt{oo4},
\texttt{oo6}, or 
\texttt{oo8})
or $b'=b+1$ (otherwise).
We have $a\geq c+1$ and in each case
$a'=a$ and $c'=c$ so $a'\geq c'+1$, thus
$h'$ is 0 or in $LL$ or in $S_1$
and we are done.

\begin{figure}[h]{\centering
\input{r3}

}\caption{How Left plays on $\mS_0$: rule 3.}
\label{fr3}
\end{figure}

Subcase 3: Rules 1 and 2 do not apply, Rule 3 applies.
Thus $y=z=0$ and $a\geq c$ and $c\geq 1$.
In case 3a, $h'=0$.
In each other case, $a'=a$ and $c'=c-1$, so $a'\geq c'+1$ and we are done.

Subcase 4: Rules 1,2,3 do not apply, Rule 4 applies.
Thus $c=y=z=0$ and $a \geq c$.
In case 4b, $h'$ in $LL$.
In each other case, $a' = a + 1$ and $c' = c$, so $a'\geq c' + 1$ thus $h'$ is in $S_1$ and we are done.

Subcase 5: 
Rules 1,\ldots\ 4 do not apply, Rule 5 applies.
Thus $b=c=y=z=0$ and $a \geq c$.
If $a=0$, In case 5g and 5j, we have $d' = 0$. Since there will be no other parts remains other than in \texttt{oo8} or \texttt{xxo}, $h'$ will either be in 0 or $LL$ or in $S_2$(if there is \texttt{xxo}).
In all other case, we have either $d' = 0$ or $d'+e' \geq 3$ and $e\geq 1$ thus $h'$ is in $S_2$.
If $a\geq 1$, since $c' = c = 0$ and $a' = a$, so $a' \geq c' + 1$ thus $h'$ is in $S_1$ and we are done.

Subcase 6: Rules 1,\ldots\ 5 do not apply, 
Rule 6 applies.
Thus $b=c=d=y=z=0$ and $a \geq 1$.
In case 6d, $a' = a - 1$ and $b=c=d=y=z=0$ and $d' = d = 0$ and $e' = e + 1 \geq 1$ thus $h'$ is in $S_2$.
In each other case, $a' = a = 1$ and $c' = c = 0$ so $a'\geq c' + 1$ thus $h'$ is in $S_1$ and we are done.

Subcase 7: Rules 1,\ldots\ 6 do not apply, 
Rule 7 applies.
Thus $a=b=c=d=y=z=0$.
In case 7a, $e' = e + 2$ and $d' = d + 1 = 1$, thus $h'$ is in $S_2$.
In case 7b, $e' = e - 1$ and $a'=b'=c'=d'= f' =y'=z'=0$ thus $h'$ is in $S_2$ or 0 if $e = 0$.
\end{proof}

\begin{figure}[h]{\centering
\input{r4}

}\caption{How Left plays on $\mS_0$: Rule 4.}
\label{fr4}
\end{figure}

\begin{figure}[h]{\centering
\input{r5}

}\caption{How Left plays on $\mS_0$: Rule 5.}
\label{fr5}
\end{figure}

\begin{figure}[h]{\centering
\input{r67}

}\caption{How Left plays on $\mS_0$: Rules 6 and 7.}
\label{fr67}
\end{figure}

\begin{proof}[Proof of Corollary.]
Let game $g$ be in $\mA$.
If $g$ is in \{2, 4, 12 \}
then
$L$ plays to 0, \myup, or 0 respectively
(2 is $*$, 12 is $4*$).
If $g$ is in \{8, 10, 14, \ldots\}
then
$L$ plays $g=2t$ to $2t-3$, which is in $\mO'$. 
In each case $L$ wins.
\end{proof}

\section*{Verification}
We implemented a Python program 
that verifies the correctness
of our method on parts with size up to 100:
for each ALC game starting from
a8, a10, \ldots, a100, Left can win.
See Figure~\ref{fver}.
The code is available on github \cite{tunglamgithub}.
Our program performs a tree search: 
when it is Right's turn, we consider all possible moves;
when it is Left's turn, we move according to Rules 1 through 7.
At each point in the search, the current game
is stored as a set of parts.
Once solved, a tree node result 
(win for nodes with Left-to-play, lose for nodes with Right-to-play)
is stored in a Python dictionary. At each node in the search,
before solving that node's game,
we check the dictionary to see whether it has already been solved.

\begin{figure}[ht!]
{\centering \includegraphics[width=.5\textwidth]{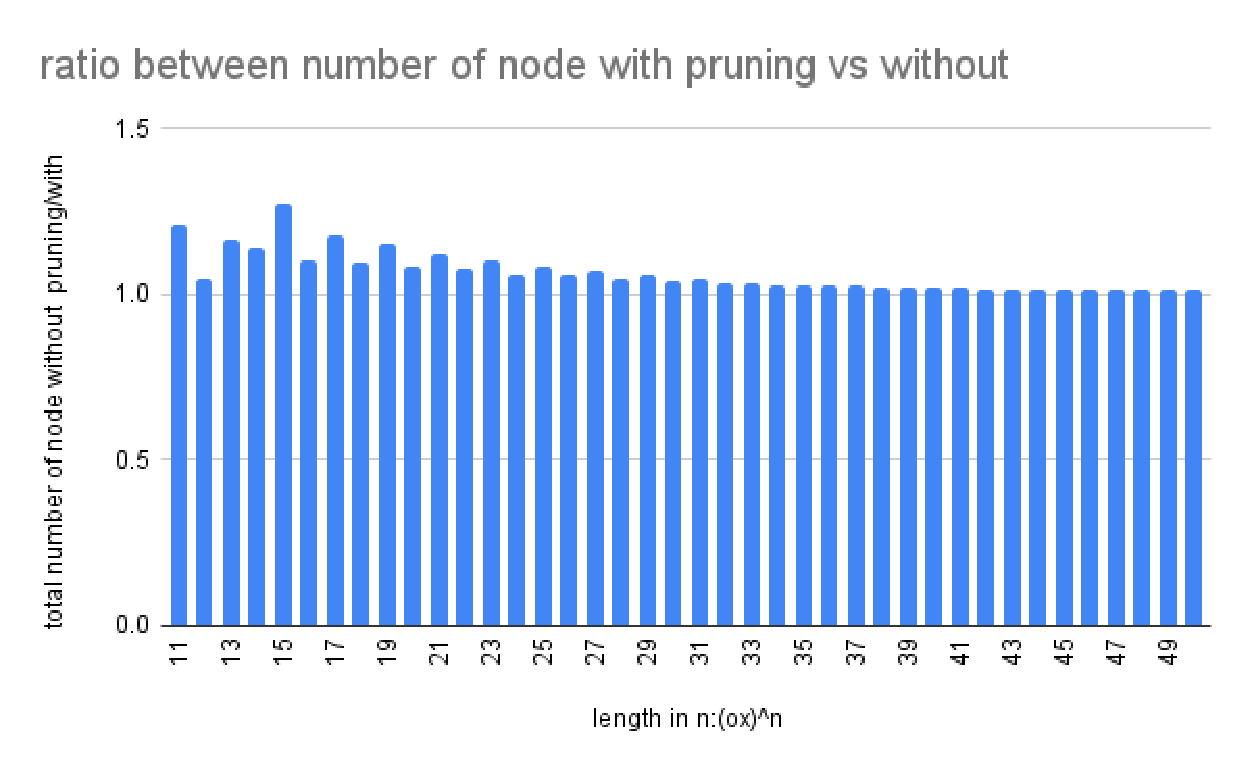}

}\caption{Node count ratio: original/improved versions.}
\label{node_ratio}
\end{figure}

\begin{figure}[ht!]
{\centering \includegraphics[width=.5\textwidth]{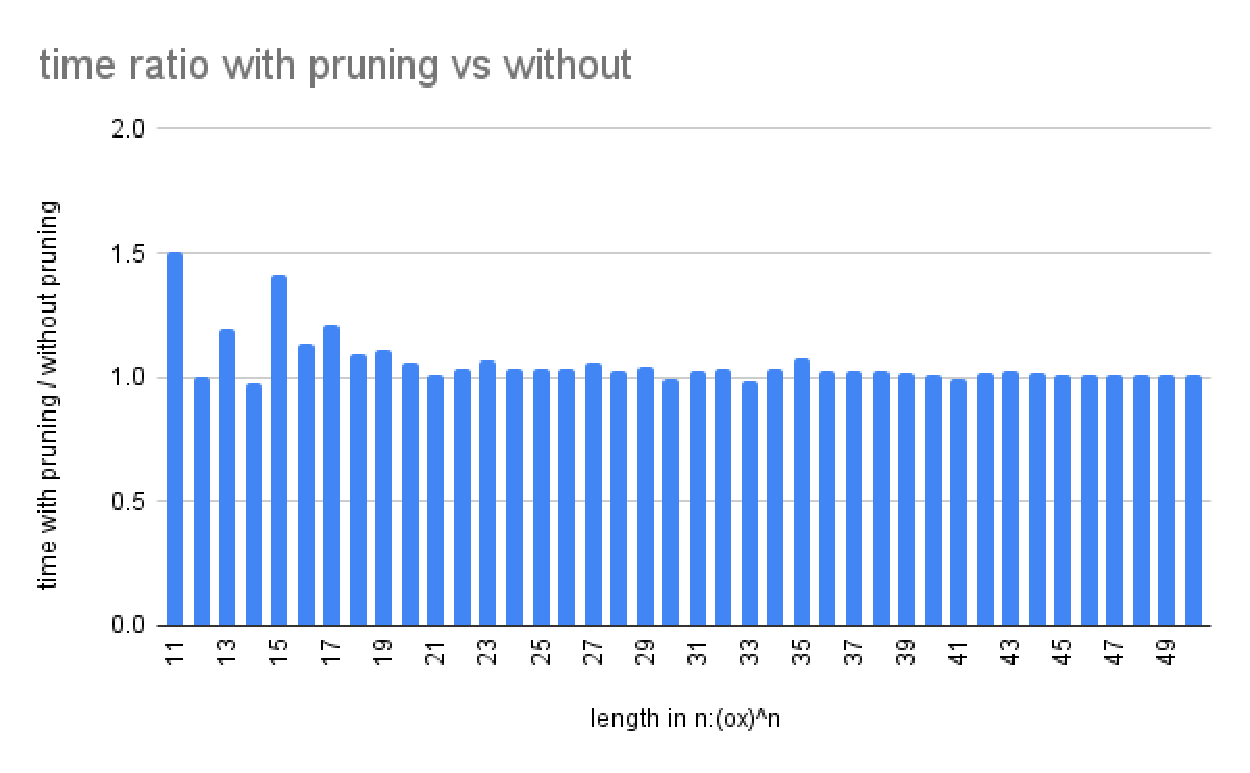}

}\caption{Time ratio: original/improved versions.}
\label{time_ratio}
\end{figure}

\begin{figure}[ht!]
{\centering \includegraphics[width=.5\textwidth]{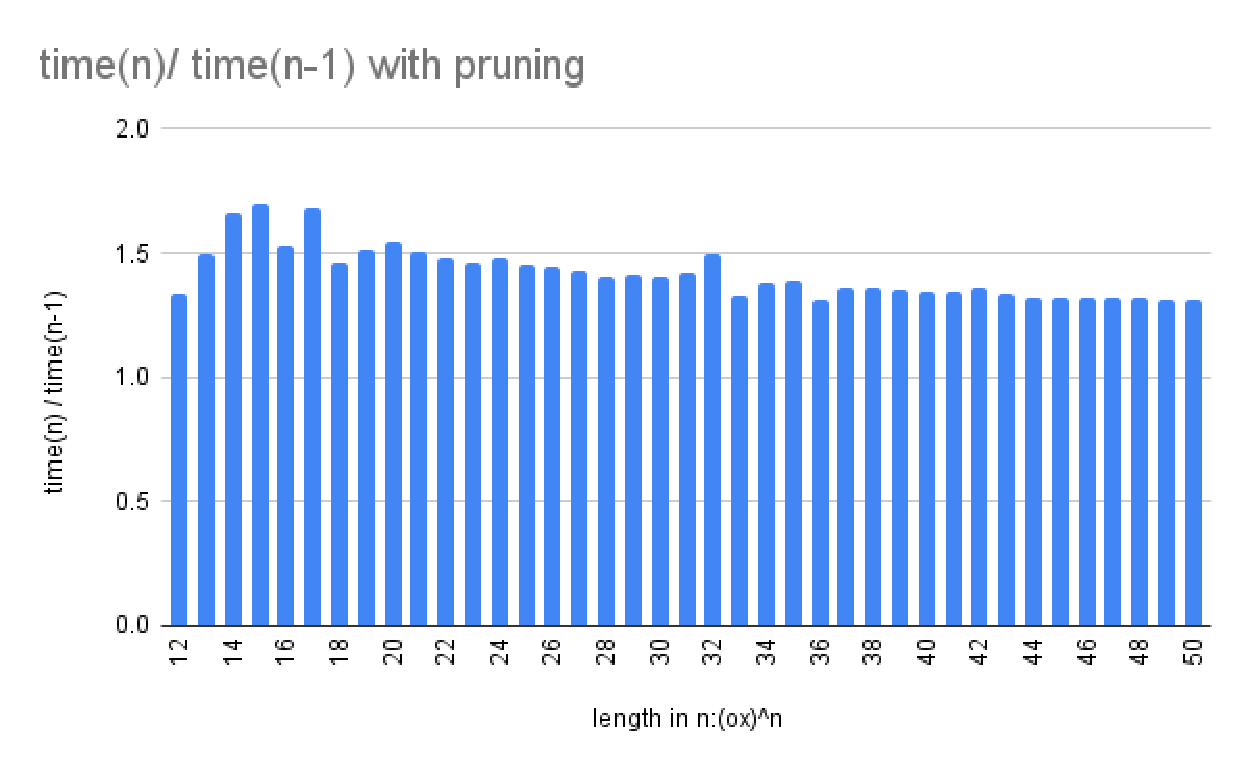}

}\caption{Runtime, improved version: time(n)/time(n-1).}
\label{time_complexity}
\end{figure}

We also implemented a program in which Left follows an improved
rule set: Left plays a(2$n$) to o(2$n-3$);
Right then plays o(2$n-3$) to  a(2$j$)$\myadd$oo(2$k$) for some $j$ and $k$;
when a(2$j$) is long enough, instead of playing a(2$j$) to o(2$j-3$),
Left plays a(2$j$) to o(2$(j-k)-1$)$\myadd$xx(2$k$),
leaving o(2$(j-k)-1$)$\myadd$xx(2$k$)$\myadd$oo(2$k$) $=$
o(2$(j-k)-1$), namely a game with just one part.
We verified the correctness
of this improved method on parts with size up to 120.
See Figures~\ref{fver2} and \ref{fver2b}.
Surprisingly,
the improved rule set reduces runtime and node count
by a negligible factor (about 1\% on a100).
See Figures~\ref{node_ratio} through \ref{time_complexity}.
From the last figure,
the runtime of our algorithm on parts of size $n$
looks to be bounded by a constant times $1.5^{n}$.

\begin{figure}[ht!]
{\tiny\centering
\input{data/ver}

}\caption{Runtime and dictionary size to verify game a(2n).}
\label{fver}
\end{figure}

\begin{figure}[ht!]
{\tiny\centering
\input{data/ver2}

}\caption{Improved ruleset, runtime and dictionary size to verify a(2n).}
\label{fver2}
\end{figure}

\begin{figure}[ht!]
{\tiny\centering
\input{data/ver2b}

}\caption{Improved ruleset, data for a102 to a120.}
\label{fver2b}
\end{figure}

\section*{Open problems}
Here are some open ALC problems:

\begin{enumerate}
\item which Left-Right-symmetric ALC games are in $\mP$?
\item is there a polytime algorithm to find the winner of an ALC game?
\end{enumerate}

\section*{Acknowledgements}
We presented this work
on January 31 2025 at
Combinatorial Game Theory Colloquium V in Lisbon.
We thank Alda Carvalho,
Carlos Pereira dos Santos,
Jorge Nuno Silva and Tiago Hirth
for organizing the conference.
We thank Martin M\"{u}ller for helpful comments
and Denilson Barbosa for providing 2025 summer research support.
\bibliographystyle{plain}
\bibliography{alc}
\end{document}

%% file: r1.tex
\begin{tabular}{|c|c|c|} \hline
rule & \{game\} or subgame & result \\ \hline
1&\multicolumn{2}{c|}{$\mA'$ non-empty}\\ \hline
a  & \{ a8 \myadd\ a2\}         
    & xxo \myadd\ a4 \myadd\ a2  (avoid  o5 \myadd\ a2) \\ \hline
b  & \{a10 \myadd\ a4\}         
    & o5 \myadd\ xxox \myadd\ a4 (avoid  o7 \myadd\ a4) \\ \hline
c  & \{a18 \myadd\ a4 \myadd\ a2\}       
    & a14 \myadd\ xxo \myadd\ a4 \myadd\ a2  (avoid o15 \myadd\ a4 \myadd\ a2) \\ \hline
d  & any of a8, a10, a14, \ldots & o5, o7, o11, \ldots \\ \hline
\end{tabular}

%% file: r2.tex
\begin{tabular}{|c|c|c|} \hline
rule & \{game\} or subgame & result \\ \hline
2&\multicolumn{2}{c|}{$\mA'$ empty, $\moA'$ non-empty}\\ \hline
  & oo7, oo9, \ldots & oo4, oo6, \ldots \\ \hline
\end{tabular}

%% file: r3.tex
\begin{tabular}{|c|c|c|} \hline
rule & \{game\} or subgame & result \\ \hline
3&\multicolumn{2}{c|}{$\mA', \moA'$ empty, $\moOo'$ non-empty}\\ \hline
a  &\{o5 $\myadd$ oox\}     & xxo $\myadd$ oox \\ \hline
b  & o5 $\myadd$ a4 $\myadd$  a2 $\myadd$  oox  & o5 $\myadd$ a4 $\myadd$ a2 $\myadd$ a2 \\ \hline
c  & a4 $\myadd$ oox        & xxo $\myadd$ oox \\ \hline
d  & oox          & a2 \\ \hline
e  & oo9oo, oo11oo, \ldots & oo4, oo6, \ldots \\ \hline
\end{tabular}

%% file: r4.tex
\begin{tabular}{|c|c|c|} \hline
rule & \{game\} or subgame & result \\ \hline
4&\multicolumn{2}{c|}{ $\mA', \moA', \moOo'$ empty, $\moO'$ non-empty}\\ \hline
a  & \{oo10 $\myadd$ a2\} & o7 $\myadd$ a2 $\myadd$ a2 
(avoid o5 $\myadd$ a2) \\ \hline
b  & \{oo12 $\myadd$ a4\} & oo8 $\myadd$ xxo $\myadd$ a4 
(avoid o7 $\myadd$ a4)\\ \hline
c  & oo10       &  o5 \\ \hline
d  & oo12       &  o7 \\ \hline
e  & oo14       &  o11 $\myadd$ a2 \\ \hline
f  & oo16       &  o11 \\ \hline
g  & oo18       &  o13 (avoid o15 $\myadd$ a4 $\myadd$ a2)\\ \hline
h  & oo20       &  o17 a2 (avoid o15 $\myadd$ a4 $\myadd$ a2)\\ \hline
i  & oo22, oo24, \ldots &  o17, o19, \ldots \\ \hline
\end{tabular}

%% file: r5.tex
\begin{tabular}{|c|c|c|} \hline
rule & \{game\} or subgame & result \\ \hline
5&\multicolumn{2}{c|}{$\mA'$, $\moA'$, $\moOo'$, $\moO'$ empty, $\mI$ non-empty}\\ \hline
a  & oo6 $\myadd$ oo6 $\myadd$ a4 & xxo $\myadd$ oo6 $\myadd$ a4\\ \hline
b  & oo6 $\myadd$ oo6 $\myadd$ a2 & xxo $\myadd$ oo6 $\myadd$ a2\\ \hline
c  & oo6 $\myadd$ oo6    & ooxo $\myadd$ oo6 $=$ xxo $\myadd$ oo6 $\myadd$ a2 \\ \hline
d  & oo6 $\myadd$ a4 $\myadd$ a2  & xxo $\myadd$ a4 $\myadd$ a2 \\ \hline
e  & oo6 $\myadd$ a4     & ooxo $\myadd$ a4 $=$ xxo $\myadd$ a4 $\myadd$ a2 \\ \hline
f  & oo6 $\myadd$ a2 & ooxo $\myadd$ a2 $=$ xxo \\ \hline
g  & a4 $\myadd$ a2 & a2 $\myadd$ a2\\ \hline
h  & oo6 & xxo \\ \hline
i  & a4  & xxo \\ \hline
j  & a2  &  0 \\ \hline
\end{tabular}

%% file: r67.tex
\begin{tabular}{|c|c|c|} \hline
rule & \{game\} or subgame & result \\ \hline
6&\multicolumn{2}{c|}{$\mA'$, $\moA'$, $\moOo'$, $\moO'$, $\mI$ empty}\\ \hline
a  & o13, o15, \ldots & o11, o13 \ldots \\ \hline
b  & o11 & o7 xxo \\ \hline
c  & o7 & o5 \\ \hline
d  & o5 & xxo \\ \hline
7&\multicolumn{2}{c|}{$\mA'$, $\moA'$, $\moOo'$, $\moO'$, $\mI$, $\mO'$ empty}\\ \hline
a  & oo8 & ooxo xxo  \\ \hline
b  & xxo & 0 \\ \hline
\end{tabular}

%% file: data/ver.tex
\begin{tabular}{|c|c|c|c|} \hline 
n & runtime in seconds & number of left node & number of right node\\ \hline
4 & 0.0 & 9 & 5\\ \hline
5 & 0.0 & 11 & 5\\ \hline
6 & 0.0 & 7 & 3\\ \hline
7 & 0.0 & 16 & 8\\ \hline
8 & 0.0 & 22 & 9\\ \hline
9 & 0.00 & 34 & 17\\ \hline
10 & 0.0 & 37 & 17\\ \hline
11 & 0.00 & 70 & 29\\ \hline
12 & 0.00 & 101 & 41\\ \hline
13 & 0.00 & 140 & 54\\ \hline
14 & 0.00 & 244 & 93\\ \hline
15 & 0.00 & 328 & 121\\ \hline
16 & 0.01 & 484 & 167\\ \hline
17 & 0.02 & 738 & 243\\ \hline
18 & 0.03 & 1050 & 328\\ \hline
19 & 0.04 & 1483 & 446\\ \hline
20 & 0.07 & 2131 & 618\\ \hline
21 & 0.11 & 2974 & 831\\ \hline
22 & 0.16 & 4144 & 1128\\ \hline
23 & 0.24 & 5732 & 1519\\ \hline
24 & 0.36 & 7889 & 2035\\ \hline
25 & 0.52 & 10773 & 2714\\ \hline
26 & 0.75 & 14734 & 3631\\ \hline
27 & 1.08 & 19846 & 4788\\ \hline
28 & 1.52 & 26691 & 6316\\ \hline
29 & 2.15 & 35784 & 8333\\ \hline
30 & 3.02 & 47670 & 10900\\ \hline
31 & 4.29 & 63220 & 14231\\ \hline
32 & 6.41 & 83717 & 18574\\ \hline
33 & 8.53 & 110323 & 24132\\ \hline
34 & 11.76 & 144847 & 31236\\ \hline
35 & 16.29 & 189611 & 40376\\ \hline
36 & 21.38 & 247575 & 52033\\ \hline
37 & 29.10 & 322116 & 66897\\ \hline
38 & 39.60 & 418243 & 85855\\ \hline
39 & 53.47 & 541192 & 109822\\ \hline
40 & 71.90 & 698675 & 140171\\ \hline
41 & 96.71 & 899744 & 178591\\ \hline
42 & 132.08 & 1155855 & 226969\\ \hline
43 & 176.71 & 1481255 & 287881\\ \hline
44 & 233.94 & 1894471 & 364467\\ \hline
45 & 308.83 & 2417306 & 460442\\ \hline
46 & 408.71 & 3078103 & 580570\\ \hline
47 & 540.34 & 3911499 & 730786\\ \hline
48 & 713.29 & 4961035 & 918104\\ \hline
49 & 936.20 & 6279679 & 1151475\\ \hline
50 & 1230.69 & 7934795 & 1441804\\ \hline
\end{tabular}

%% file: data/ver2.tex
\begin{tabular}{|c|c|c|c|} \hline 
n & runtime in seconds & number of left node & number of right node\\ \hline
4 & 0.0 & 9 & 5\\ \hline
5 & 0.0 & 11 & 5\\ \hline
6 & 0.0 & 7 & 3\\ \hline
7 & 0.00 & 16 & 8\\ \hline
8 & 0.0 & 22 & 9\\ \hline
9 & 0.0 & 34 & 17\\ \hline
10 & 0.00 & 36 & 17\\ \hline
11 & 0.00 & 57 & 25\\ \hline
12 & 0.00 & 96 & 40\\ \hline
13 & 0.00 & 119 & 47\\ \hline
14 & 0.00 & 213 & 83\\ \hline
15 & 0.00 & 256 & 96\\ \hline
16 & 0.01 & 436 & 155\\ \hline
17 & 0.01 & 617 & 216\\ \hline
18 & 0.02 & 953 & 304\\ \hline
19 & 0.04 & 1275 & 401\\ \hline
20 & 0.07 & 1957 & 581\\ \hline
21 & 0.11 & 2635 & 763\\ \hline
22 & 0.16 & 3842 & 1061\\ \hline
23 & 0.22 & 5167 & 1402\\ \hline
24 & 0.34 & 7430 & 1942\\ \hline
25 & 0.50 & 9896 & 2546\\ \hline
26 & 0.73 & 13928 & 3456\\ \hline
27 & 1.02 & 18421 & 4508\\ \hline
28 & 1.48 & 25445 & 6054\\ \hline
29 & 2.07 & 33670 & 7917\\ \hline
30 & 3.04 & 45820 & 10522\\ \hline
31 & 4.18 & 60439 & 13702\\ \hline
32 & 6.18 & 80944 & 18003\\ \hline
33 & 8.65 & 106385 & 23372\\ \hline
34 & 11.36 & 140699 & 30379\\ \hline
35 & 15.15 & 183948 & 39265\\ \hline
36 & 20.78 & 241464 & 50759\\ \hline
37 & 28.42 & 313963 & 65259\\ \hline
38 & 38.68 & 409244 & 83939\\ \hline
39 & 52.31 & 529291 & 107362\\ \hline
40 & 71.01 & 685317 & 137291\\ \hline
41 & 97.02 & 882347 & 174916\\ \hline
42 & 129.72 & 1135756 & 222573\\ \hline
43 & 172.50 & 1455253 & 282308\\ \hline
44 & 230.06 & 1864274 & 357832\\ \hline
45 & 305.91 & 2378531 & 452072\\ \hline
46 & 405.76 & 3031968 & 570444\\ \hline
47 & 536.10 & 3853219 & 718202\\ \hline
48 & 707.37 & 4891098 & 902851\\ \hline
49 & 930.19 & 6191146 & 1132524\\ \hline
50 & 1219.64 & 7828532 & 1418867\\ \hline
\end{tabular}

%% file: data/ver2b.tex
\begin{tabular}{|c|c|c|c|} \hline 
n & runtime in seconds & number of left node & number of right node\\ \hline
51 & 1597.32 & 9873400 & 1773823\\ \hline
52 & 2079.01 & 12437109 & 2214530\\ \hline
53 & 2719.68 & 15632837 & 2759922\\ \hline
54 & 3530.88 & 19623639 & 3434831\\ \hline
55 & 4581.00 & 24584029 & 4267453\\ \hline
56 & 5926.90 & 30759481 & 5295222\\ \hline
57 & 8274.98 & 38415118 & 6559902\\ \hline
58 & 11022.62 & 47913880 & 8116178\\ \hline
59 & 13444.17 & 59661929 & 10027059\\ \hline
60 & 18485.15 & 74194971 & 12372128\\ \hline
\end{tabular}